\documentclass[10pt, twocolumn, report]{amsart} 
\usepackage{amsmath,amscd, amsthm,  amssymb, bm, pifont }

\date{}
\newtheorem{thm}{Theorem}[section]
\newtheorem{prop}[thm]{Proposition}
\newtheorem{lem}[thm]{Lemma}
\newtheorem{rem}[thm]{Remark}

\newtheorem{prob}[thm]{Problem}  
\numberwithin{equation}{section}

\def\Im{\operatorname{Im}}
\newcommand{\Z}{\mathbb{Z}}
\newcommand{\Q}{\mathbb{Q}}
\newcommand{\C}{\mathbb{C}}

\begin{document}
\title{ The integral cohomology ring of $E_{8}/T$  }
\author{Masaki Nakagawa}
\pagestyle{plain}



\thanks{
$^*$ Partially supported by the Grant-in-Aid for Scientific  Research 
(C)  21540104,   Japan Society for the Promotion of Science. 
}

\address{Department of General Education \endgraf
                      Kagawa National College of Technology \endgraf
                      Takamatsu 761-8058 \\ Japan}
\email{nakagawa@t.kagawa-nct.ac.jp}

\twocolumn[
\maketitle
]

\begin{abstract}
     We give a complete description of the integral cohomology ring of the flag 
    manifold $E_8/T$, where $E_8$ denotes the compact exceptional Lie group of rank $8$ and $T$ 
    its maximal torus,  by the method due to Borel and Toda.  This completes the computation of the 
    integral cohomology rings of the flag manifolds for all compact connected simple Lie groups.  
\end{abstract} 

\section{Introduction} 
Let $G$ be a compact connected Lie group and $T$ a maximal torus of $G$.  Then the homogeneous space
$G/T$ is called a ({\it full} or {\it complete}) {\it  flag manifold} and plays an important role in modern mathematics.   
In algebraic topology, the following problem is classical:  
\begin{prob}  Determine the integral cohomology ring of the flag manifold $G/T$ for $G$ 
a compact  connected simple Lie group. 
\end{prob} 
The computation of the integral cohomology ring of $G/T$ was started by  Borel in 1953 (\cite{Bor53}).  
Borel considered the spectral sequence for the  fibration 
\[    
      G/T \overset{\iota}{\longrightarrow}  BT \overset{\rho}{\longrightarrow}  BG,    
\]
where $BT$   (resp. $BG$) denotes the classifying space of $T$ (resp. $G$), 
and obtained the following description of the cohomology ring of $G/T$; Let $k$ be a field of characteristic $p$, and  
suppose that $p = 0$ or $G$ has no $p$-torsion.  The Weyl group $W$ of $G$ acts naturally on $T$, and hence, 
on $BT$ and also on $H^{*}(BT;k)$.    Then the homomorphism 
\[ 
     \iota^{*}: H^{*}(BT;k) \longrightarrow H^{*}(G/T;k)
\]  induces   the following  isomorphism:
\begin{equation}  \label{eqn:Borel}     
      H^*(BT;k)/(H^{+}(BT;k)^{W})   \longrightarrow      H^*(G/T;k),   
\end{equation} 
where $(H^{+}(BT;k)^{W})$ denotes the ideal of $H^*(BT;k)$ generated by  $W$-invariant polynomials  of 
positive degrees.   This is the {\it Borel presentation} of the cohomology ring of $G/T$.  
Borel's  method is valid for the integer coefficients when $G$ has no torsion.  
So  the cases of $SU(n)$ and $Sp(n)$ follow immediately. However, when $G$ has $p$-torsion, Borel's result 
does not hold.  In 1955, Bott and Samelson developed 
an algorithm  for computing the integral cohomology ring of $G/T$ by means of  the so-called ``Bott-Samelson $K$-cycle''\footnote{A  Bott-Samelson $K$-cycle is refered to as a {\it Bott tower} in toric topology, and a  {\it Bott-Samelson-Demazure-Hansen variety} in 
representation theory.},  and determined 
the case of the exceptional group $G_{2}$ explicitly (\cite{Bot-Sam55}, \cite[Theorem  III']{Bot-Sam58}). 
The above problem could  be solved theoretically  by their method.  However, it seems difficult to apply 
this  to  other exceptional goups $F_{4}$, $E_{6}$, $E_{7}$ and $E_{8}$.  In 1975,   Toda gave another 
useful  description of  the integral cohomology ring of $G/T$ from the known results on  the mod $p$ cohomology 
rings $H^*(G;\Z/p\Z)$ of $G$ for all primes $p$ and the rational cohomology ring $H^*(G/T;\Q)$ of $G/T$ 
(\cite[Theorem 2.1, Proposition 3.2]{Toda75}).
Based on  Toda's method,  the cases of $SO(2n)$, $SO(2n+1)$  were settled by  Toda and  Watanabe 
(\cite[Theorem 2.1, Corollary 2.2]{Toda-Wat74}).\footnote{The results for $SO(2n)$, $SO(2n+1)$  also appeared in \cite{Mar74} by computing the ``cohomology 
ring of the root system'' due to   Demazure  \cite{Dem73}.}    The cases of $F_{4}$ and $E_{6}$
were also settled  in \cite[Theorems A, B]{Toda-Wat74}. The case of $E_{7}$ was settled  by the author (\cite[Theorem 5.9]{Nak01}). 
So the only remaining case is $G = E_{8}$.  In this article, we determine the integral cohomology ring of $E_8/T$ explicitly 
along the line of Toda's method.\footnote{Recently Duan and Zhao also computed it in the context of the {\it Schubert calculus} (\cite{DZ08}).}  
This completes the computation of the integral 
cohomology rings of the flag manifolds for all compact connected simple Lie groups.

\section{Ring of invariants of $W(E_{8})$}
Let $E_{8}$ be the compact simply connected simple exceptional Lie group of rank $8$
and $T$ a maximal torus.  Following  \cite{Bou68}, we take the simple roots $\{ \alpha_{i} \}_{1 \leq i \leq 8}$
and denote by  $\{ \omega_{i} \}_{1 \leq i \leq 8}$ the  corresponding fundamental weights. 
In topology,  it is customary that  roots and weights are regarded as elements of $H^{2}(BT;\Z)$. 
Let $s_{i} \; (1 \leq i \leq 8)$ denote the simple reflection corresponding to the simple root $\alpha_{i}$.
Then the Weyl group $W(E_{8})$ of $E_{8}$ is generated by simple reflections $s_{i} \; (1 \leq i \leq 8)$. 
As in   \cite[\S 2]{Nak09}, we put 
\begin{equation}  \label{eqn:t_i}
     \begin{array}{lll} 
          t_{8} & = \omega_{8},  \; t_{i} = s_{i+1}({t_{i+1}}) \; (2 \leq i \leq 7),  \medskip  \\ 
          t_{1} & = s_{1}(t_{2}), \;   t      = \omega_{2}.  \medskip  
     \end{array}    
\end{equation} 
Then we have 
\[  
        H^*(BT;\Z) = \Z[t_{1}, \ldots, t_{8}, t]/(c_{1} - 3t)    
\]
for $c_{i} = e(t_{1}, \ldots, t_{8})$, the  $i$-th elementary symmetric polynomial in $t_{1}, \ldots, t_{8}$.  

According to Chevalley  \cite{Che55}, the ring of invariants of the Weyl group $W(E_{8})$  over $\Q$  is generated by 
$8$ algebraically   independent polynomials ({\it basic invariants}) of degrees $2, 8, 12, 14, 18, 20, 24, 30$.  
By computing the Chern character of the adjoint representation of $E_{8}$, of dimension $248$, 
we obtain  the basic invariants  $I_{j} \; (j = 2, 8, 12, 14, 18, 20, 24, 30)$  explicitly  in  \cite[Lemma 2.3]{Nak09}
(see also \cite[2.3]{Meh88}).   Thus we have the following:  
\begin{lem}  \label{lem:W(E_8)-invariants}  
The ring of invariants of the Weyl group $W(E_{8})$ over $\Q$ is given by 
\[ 
    H^*(BT;\Q)^{W(E_{8})} = \Q[I_{2}, I_{8}, I_{12}, I_{14}, I_{18}, I_{20}, I_{24}, I_{30}].  
\]  
\end{lem} 
By (\ref{eqn:Borel}) and Lemma \ref{lem:W(E_8)-invariants}, we can compute the rational 
cohomology ring $H^*(E_{8}/T;\Q)$  of $E_{8}/T$.

\section{Integral cohomology ring  of $E_8/T$}
As mentioned in the introduction, we compute \\  $H^*(E_8/T;\Z)$ following Toda's method. 
Since $E_8$ is simply connected, 
the homomorphism 
\[  \iota^*: H^{2}(BT;\Z)    \longrightarrow     H^{2}(E_{8}/T;\Z)
\] is an isomorphism.   
Under this isomorphism, we denote the $\iota^*$-images of $t_{i}$ and $t$ by the same symbols. 
Thus $H^2(E_8/T;\Z)$ is a free $\Z$-module generated by $t_{i} \; (1 \leq i \leq 8)$ and $t$
with a relation $c_{1} = 3t$. 
 Toda \cite{Toda75} gave a   general description of the integral cohomology ring 
 of $E_8/T$.  In our situation, his result is stated as follows:
 \begin{prop} [\cite{Toda75}, Proposition 3.2] \label{prop:E_8/T}
 The integral cohomology ring of $E_{8}/T$ is of the form$:$  
 \begin{align*}
    &  H^{*}(E_{8}/T;\Z)  \\
           = & \; \dfrac{\Z[t_{1},\ldots,t_{8}, t, \gamma_{3},  \gamma_{4}, \gamma_{5}, \gamma_{6}, \gamma_{9},
                \gamma_{10},  \gamma_{15}]}
             {\left (  \begin{array}{lll} 
                            &  \rho_{1},  \rho_{2},  \rho_{3},  \rho_{4}, \rho_{5},  \rho_{6},  \rho_{8},  \rho_{9},
                               \rho_{10}, \rho_{12}, \rho_{14}, \rho_{15},   \\
                            & \rho_{18},\rho_{20},  \rho_{24},\rho_{30}
                       \end{array}   
              \right ),  }
   \end{align*}
where  $t_{1}, \ldots, t_{8}, t \in H^{2}$ are as above, and  $\gamma_{i} \in H^{2i} \; (i = 3,4,5,6,9,10,15)$, 
and 
 \begin{align*}
  \rho_{1} &= c_{1} - 3t, \\
  \rho_{i} &= \delta_{i} - 2\gamma_{i} \; (i= 3,5,9,15),  \\
  \rho_{i} &= \delta_{i}-3\gamma_{i}\; (i = 4,10), \\ 
  \rho_{6} &= \delta_{6}-5\gamma_{6}.   
  \end{align*}
Here $\delta_{i} \; (i = 3,4,5,6,9,10,15)$ are  arbitrary elements of $H^*(E_{8}/T;\Z)$  satisfying 
  \begin{align*}
 \delta_{3} &\equiv Sq^{2}(\rho_{2}) ,\; \delta_{5} \equiv Sq^{4}(\delta_{3}),
  \; \delta_{9} \equiv Sq^{8}(\delta_{5}), \\
 \delta_{15} & \equiv Sq^{14}(\rho_{8})   \mod 2,  \\
 \delta_{4} &\equiv  \mathcal{P}^{1}(\rho_{2}),\; \delta_{10} \equiv 
   \mathcal{P}^{3}(\delta_{4})   \mod 3,  \\
 \delta_{6} &\equiv \mathcal{P}^{1}(\rho_{2}) \mod 5.  
  \end{align*}
Other relations $\rho_{j} \; (j = 2,8,12,14,18,20,24,30)$ are  determined by the 
maximum of the integers $n_{j}$ in 
 \begin{equation} \label{eqn:rho_j}
    n_{j}  \cdot \rho_{j} \equiv \iota^{*}(I_{j})   \mod (\rho_{i}; i < j),          
 \end{equation}
where $I_{j} \; (j = 2, 8, 12, 14, 18, 20, 24, 30)$ are the  basic  invariants 
of the Weyl group $W(E_{8})$ given in Lemma $\ref{lem:W(E_8)-invariants}$. 
\end{prop} 
We will carry out his program for $E_{8}$.  Fortunately,  in  \cite[Lemma 4.2]{Nak09}, partial computation of 
$H^*(E_8/T;\Z)$ has been  made up to degrees $36$. So  we need only to determine the higher relations 
$\rho_{20}$, $\rho_{24}$ and $\rho_{30}$ explicitly.   However, it seems  difficult to compute them 
directly from the basic invariants $I_{20}$, $I_{24}$ and $I_{30}$.  So we 
make use of a certain subgroup of $E_8$. Namely,  let $C$ be the centralizer of a  one dimensional 
torus determined by $\alpha_{i} = 0 \; (i \neq 8)$.  The local  type of $C$  can be read off from 
the Dynkin diagram of $E_{8}$ (\cite{Bor-Sie49}), and we have, in fact, 
  \[  C = T^1 \! \cdot \! E_7 \quad \text{and} \quad T^1 \cap E_{7} \cong  \Z/2\Z,   \]
where $T^{1}$ denotes a certain  one-dimensional torus. 
Consider the fibration
\[     
        E_{7}/T' \cong C/T \overset{i}{\longrightarrow} E_{8}/T
         \overset{p}{\longrightarrow}E_{8}/C,   
\]
where $T'$ is a maximal torus of $E_7$.   Since $H^{*}(E_{8}/C;\Z)$ and $H^{*}(E_{7}/T';\Z)$ have no torsion and 
vanishing odd dimensional part by  Bott \cite[Theorem A]{Bot56}, the Serre spectral sequence with integer coefficients 
for the above fibration  collapses,  and the following  sequence 
\begin{align*} 
    \Z \longrightarrow  & H^{*}(E_{8}/C;\Z) \overset{p^{*}}
   {\longrightarrow} H^{*}(E_{8}/T;\Z)   \\
    \overset{i^{*}}{\longrightarrow}
   & H^{*}(C/T;\Z) \cong H^{*}(E_{7}/T';\Z) \longrightarrow \Z  
\end{align*}   
is co-exact. In particular, $p^*$ is a split monomorphism so that $H^*(E_8/T;\Z)$ contains 
$H^*(E_8/C;\Z)$ as a direct summand if we identify $\Im p^*$ with $H^*(E_8/C;\Z)$.   

The integral cohomology ring of $E_8/C$ is determined in  \cite{Nak09}, which we 
now recall:  
\begin{thm}[\cite{Nak09}, Theorem 4.7]  \label{thm:E_8/C}
   The integral cohomology ring of $E_{8}/C$ is given as follows$:$
          \[   H^*(E_{8}/C;\Z) = \Z [u, v, w, x]/(r_{15}, r_{20}, r_{24}, r_{30}),  \] 
 where $\deg u = 2, \; \deg v = 12, \;  \deg w = 20, \; \deg x = 30$ and  
  \begin{align*}
     r_{15}   &= u^{15} - 2x, \\
     r_{20}   &= 9u^{20} + 45u^{14}v + 12u^{10}w + 60u^8v^2 + 30u^4vw  \\
              &+ 10u^2v^3 + 3w^2, \\
     r_{24}   &=  11u^{24} + 60u^{18}v +  21u^{14}w + 105u^{12}v^2 + 60u^8vw  \\
              &+ 60u^{6}v^3 + 9u^4 w^2 + 30 u^2 v^2 w + 5v^4,   \\
     r_{30}   &= -9x^2 - 12u^{9}vx  -6u^{5}wx  + 9u^{14}vw -10u^{12}v^3   \\
              &-3u^{10}w^2  + 30u^8 v^2w -35u^{6}v^{4} + 6u^{4}vw^{2}  \\
              &- 10u^2v^3w -4v^5 -2w^{3}.
\end{align*}   
\end{thm}   

\begin{rem} 
 The integral cohomology ring of $E_{8}/C$ is also computed by Duan and Zhao  \cite[Theorem 7]{DZ05} in the context of the 
 Schubert calculus.  In our forthcoming paper $($\cite{Kaji-Nak1}$)$, we will show that Theorem $\ref{thm:E_8/C}$  completely coincides with the result of Duan 
 and Zhao by means of the {\it divided difference operators} due to Bernstein-Gelfand-Gelfand \cite{BGG73} 
 and Demazure \cite{Dem73}. 
\end{rem} 

  The relations $r_{20}$, $r_{24}$ and $r_{30}$ of $H^*(E_8/C;\Z)$ correspond to the relations 
  $\rho_{20}$, $\rho_{24}$ and $\rho_{30}$ of $H^*(E_8/T;\Z)$ respectively.  In order to make the description of 
  $H^*(E_8/T;\Z)$ complete, we have to specify the elements $u$, $v$, $w$ and $x$ explicitly in the ring 
  $H^*(E_8/T;\Z)$.  This has been accomplished  in \cite[4.2 and 6.1]{Nak09}.  
  Furthermore, since $H^*(E_8/C;\Z)$ is a direct summand of 
  $H^*(E_8/T;\Z)$, $r_{20}$, $r_{24}$ and $r_{30}$ are not divisible by any integer in the ring 
  $H^*(E_8/T;\Z)$.  So we can replace $\rho_{20}, \rho_{24}$ and $\rho_{30}$ with $r_{20}, r_{24}$ and $r_{30}$ 
  respectively. Summing up the results so far, we obtain the following main result of this article:
\begin{thm}     \label{thm:E_8/T}  
The integral cohomology ring of $E_8/T$  is given as follows$:$
\begin{align*}
    &  H^{*}(E_{8}/T;\Z)  \\
           = & \; \dfrac{\Z[t_{1},\ldots,t_{8}, t, \gamma_{3},  \gamma_{4}, \gamma_{5}, \gamma_{6}, \gamma_{9},
                \gamma_{10},  \gamma_{15}]}
               {\left (  \begin{array}{lll} 
                               &  \rho_{1},  \rho_{2},  \rho_{3},  \rho_{4}, \rho_{5},  \rho_{6},  \rho_{8},  \rho_{9},
                                  \rho_{10}, \rho_{12}, \rho_{14}, \rho_{15},   \\
                               & \rho_{18},\rho_{20},  \rho_{24},\rho_{30}
                         \end{array}   
                \right ),  }
   \end{align*}
 where $t_{1},\ldots,  t_{8},t \in H^{2}$ are as in \S $2$, $\gamma_{i} \in H^{2i} \;
  (i = 3,4,5,6,9,10, 15)$.  The relations are  
  \begin{align*}
     \rho_{1} &= c_{1}-3t, \\
     \rho_{2} &= c_{2}-4t^{2}, \\ 
     \rho_{3} &= c_{3}-2\gamma_{3}, \\    
     \rho_{4} &= c_{4}+2t^{4}-3\gamma_{4}, \\
     \rho_{5} &= c_{5}-3t\gamma_{4} + 2t^{2}\gamma_{3}-2\gamma_{5},  \\
     \rho_{6} &= c_{6}-2\gamma_{3}^{2}-t\gamma_{5}+t^{2}\gamma_{4}-2t^{6}-5\gamma_{6},   \\
     \rho_{8} &= -3c_{8}+3\gamma_{4}^{2}-2\gamma_{3}\gamma_{5}+ t(2c_{7}-6\gamma_{3}\gamma_{4})  \\
              &+ t^2(2\gamma_{3}^{2}-5\gamma_{6}) 
               +3t^{3}\gamma_{5}+4t^{4}\gamma_{4}-6t^{5}\gamma_{3}+t^{8},  \\
     \rho_{9} &= 2c_{6}\gamma_{3}+tc_{8}+t^{2}c_{7}-3t^{3}c_{6}-2\gamma_{9}, \\
     \rho_{10} &= \gamma_{5}^{2}-2c_{7}\gamma_{3}  -t^{2}c_{8}  +3t^{3}c_{7}  -3\gamma_{10}, \\
     \rho_{12} &= 15 \gamma_{6}^2 + 2 \gamma_3 \gamma_4 \gamma_5 -2 c_{7} \gamma_{5} + 2 \gamma_{3}^4 
                  + 10 \gamma_{3}^2 \gamma_6  \\
               & - 3 c_8 \gamma_4 - 2 \gamma_{4}^3  \\
               & + t (c_8 \gamma_3 -2 \gamma_{3}^2 \gamma_5 + 4 c_7 \gamma_4 + 6 \gamma_{3} \gamma_{4}^2)  \\
               &+ t^2 (3 \gamma_{10} - 25 \gamma_4 \gamma_6 - c_7 \gamma_3 -16  \gamma_{3}^2 \gamma_4) \\
               & + t^3 (25 \gamma_3 \gamma_6 - 3 \gamma_4 \gamma_5 + 10  \gamma_{3}^3)  \\
               &  + t^4 (3c_8 + 3 \gamma_3 \gamma_5 + 5 \gamma_{4}^2)  \\
               & + t^5 (-3 c_7 - 5 \gamma_3 \gamma_4) + 4 t^6  \gamma_{3}^2 - 7 t^8 \gamma_4 + 4 t^9 \gamma_3, \\
     \rho_{14} &= c_{7}^2 - 3 c_8 \gamma_6 + 6 \gamma_4 \gamma_{10} - 4c_8 \gamma_{3}^2 + 6 c_{7} \gamma_3 \gamma_4  
                  - 6  \gamma_{3}^2 \gamma_{4}^2   \\
               &- 12  \gamma_{4}^2 \gamma_6 - 2 \gamma_3 \gamma_5 \gamma_6 \\
               &  + t(24 \gamma_3 \gamma_4 \gamma_6 - 8c_{7}  \gamma_{3}^2  - 8 c_7 \gamma_6 + 4 c_{8} \gamma_5 
                  - 6 \gamma_{3} \gamma_{10}  \\
               &+ 12  \gamma_{3}^3 \gamma_4)  \\
               &  + t^2 (-2 \gamma_3 \gamma_4 \gamma_5 + 6  \gamma_{4}^3 + 2  \gamma_{3}^2 \gamma_6 + 20  \gamma_{6}^2 
                          - 4  \gamma_{3}^4   - c_7 \gamma_5)  \\
               &  + t^3 (-12 \gamma_3  \gamma_{4}^2 + 8c_8 \gamma_3 - 5 c_7 \gamma_4 + 3\gamma_5 \gamma_6) \\
               & + t^4 (3 \gamma_{10} - 26 \gamma_4 \gamma_6 + 6c_7 \gamma_3 - 4 \gamma_{3}^2 \gamma_4) \\
               &  + t^5 (24 \gamma_3 \gamma_6 + 3 \gamma_4 \gamma_5 + 12  \gamma_{3}^3) 
                  + t^6 (-6 c_8 + 2 \gamma_{4}^2)  \\
              & - 4 t^7 c_7 + t^8 (6 \gamma_6 - 6  \gamma_{3}^2) 
                - 6 t^{10} \gamma_4 + 12 t^{11} \gamma_3 - 2t^{14}, \\
     \rho_{15} &= (c_8 - t^2 c_6 + 2 t^3 \gamma_5 + 3 t^4 \gamma_4 - t^8)(c_7 - 3tc_6)  \\
               & - 2(\gamma_{3}^2 + c_6)(\gamma_9 - c_6 \gamma_3)  - 2\gamma_{15}, \\
    \rho_{18}  &=  \gamma_{9}^2 - 9c_{8}\gamma_{10} - 6\gamma_{4}^2\gamma_{10} - 4\gamma_{3}^3\gamma_{9} 
                    - 10  \gamma_{3} \gamma_{6}\gamma_{9}  \\
               &+ 2\gamma_{3}\gamma_{5}\gamma_{10}   - 2\gamma_{3}\gamma_{4}\gamma_{5}\gamma_{6} 
                   - 6c_{7}\gamma_{3}\gamma_{4}^2 + 3c_{8}\gamma_{4}\gamma_{6}  \\
               & + c_{8}\gamma_{3}^2\gamma_{4}
                   + 6\gamma_{3}^2\gamma_{4}^3   + 12\gamma_{4}^3 \gamma_{6}+ 2c_{7}^2\gamma_{4} 
                   + 2c_{7}\gamma_{3}^2\gamma_{5}   \\
               &  - 2\gamma_{3}^3\gamma_{4}\gamma_{5}+ 2c_{7}\gamma_{5}\gamma_{6} + 4\gamma_{3}^6  
                  - 10 \gamma_{6}^3  + 18\gamma_{3}^4 \gamma_{6}  \\
               & + 15 \gamma_{3}^2 \gamma_{6}^2  - 9c_{7}c_{8}\gamma_{3} \\
              & + t(-2\gamma_{3}\gamma_{5}\gamma_{9}  - 24c_{7}\gamma_{4} \gamma_{6} + 8c_{8}\gamma_{4}\gamma_{5}
                + 4c_{7}\gamma_{3}^2\gamma_{4}    
\end{align*} 
\begin{align*} 
              & + 4c_{7}\gamma_{10}   - c_{8}\gamma_{9} + 2c_{7}^2 \gamma_{3}  + 4c_{8}\gamma_{3} \gamma_{6} 
                + 12 \gamma_{3} \gamma_{4} \gamma_{10} \\
              & - 36\gamma_{3}\gamma_{4}^2 \gamma_{6}   + 12 \gamma_{3}^2\gamma_{5} \gamma_{6} 
                + c_{8}\gamma_{3}^3 + 6\gamma_{3}^4\gamma_{5} - 18\gamma_{3}^3 \gamma_{4}^2)  \\
             &  + t^2 (24\gamma_{3}^4\gamma_{4} - 2c_{8}^2 - c_{7}\gamma_{9} - 11\gamma_{3}^2\gamma_{10} 
                   + 2\gamma_{3}\gamma_{4}\gamma_{9}   \\
             & - 2c_{8}\gamma_{3}\gamma_{5}   + 16c_{7}\gamma_{3} \gamma_{6} 
                   - 3c_{7}\gamma_{4}\gamma_{5}  
                + 75\gamma_{4} \gamma_{6}^2 - 6\gamma_{4}^4   \\
             & - 9c_{8}\gamma_{4}^2      + 81 \gamma_{3}^2\gamma_{4} \gamma_{6} - 13 \gamma_{6}\gamma_{10} 
                   + 4\gamma_{3}\gamma_{4}^2 \gamma_{5}  - c_{7}\gamma_{3}^3)  \\
               & + t^3(-3\gamma_{5}\gamma_{10} - 150\gamma_{3} \gamma_{6}^2 - 135\gamma_{3}^3 \gamma_{6}
                 + 6\gamma_{3}^2 \gamma_{9}  \\
             &  - 2 c_{7}\gamma_{3}\gamma_{5}   + 21c_{7}\gamma_{4}^2  + 15c_{7}c_{8}   + 3\gamma_{4}\gamma_{5} \gamma_{6} 
                      - 3\gamma_{3}^2 \gamma_{4} \gamma_{5}  \\
             &  + 18\gamma_{3} \gamma_{4}^3   + 15 \gamma_{6}\gamma_{9} 
                       + 14c_{8} \gamma_{3} \gamma_{4} - 30\gamma_{3}^5)   \\
               & + t^4(-13c_{8} \gamma_{6} + 2\gamma_{4}\gamma_{10} - 5c_{7}^2 - 33\gamma_{3}^2 \gamma_{4}^2 
                 + 3\gamma_{5} \gamma_{9}  \\
               &  -28 \gamma_{3}\gamma_{5} \gamma_{6} - 45\gamma_{4}^2 \gamma_{6}  
                  - 41c_{7}\gamma_{3}\gamma_{4} -13 \gamma_{3}^3 \gamma_{5} - 9c_{8} \gamma_{3}^2) \\
               & + t^5(3c_{7} \gamma_{6} - 6\gamma_{4}^2\gamma_{5} + 23c_{7}\gamma_{3}^2  
                 + 105\gamma_{3}\gamma_{4} \gamma_{6} - 6c_{8} \gamma_{5}  \\
              & - 3\gamma_{4}\gamma_{9} 
                 + 45 \gamma_{3}^3 \gamma_{4})   \\
               & + t^6(11\gamma_{4}^3 - 4\gamma_{3}\gamma_{9} + 4c_{7}\gamma_{5} + 9\gamma_{3}\gamma_{4}\gamma_{5} 
                 + 12\gamma_{3}^4   \\
               & + 66\gamma_{3}^2 \gamma_{6}  + 75\gamma_{6}^2 + 2c_{8}\gamma_{4})  \\
               & + t^7(-33\gamma_{3}\gamma_{4}^2 + 12 \gamma_{3}^2 \gamma_{5}  + 15\gamma_{5} \gamma_{6})  \\ 
               & + t^8(-4\gamma_{10} + 21\gamma_{3}^2 \gamma_{4} - 5c_{7}\gamma_{3} -3\gamma_{4} \gamma_{6}) \\
               & + t^9(6\gamma_{9}  - 42\gamma_{3}^3 - 99\gamma_{3} \gamma_{6}) \\
               & + t^{10} (-4c_{8} - 6\gamma_{4}^2 -13 \gamma_{3}\gamma_{5}) \\
               & + t^{11}(3c_{7} + 27\gamma_{3}\gamma_{4}) 
                 +t^{12} (60 \gamma_{6} + 18\gamma_{3}^2)  \\
               & + 6t^{13} \gamma_{5} - 9t^{14} \gamma_{4} 
                 - 12t^{15} \gamma_{3} +10t^{18}, \\
         \rho_{20}   &= 9u^{20} + 45u^{14}v + 12u^{10}w + 60u^8v^2 + 30u^4vw  \\
                     &  + 10u^2v^3 + 3w^2, \\
         \rho_{24}   &=  11u^{24} + 60u^{18}v +  21u^{14}w + 105u^{12}v^2 + 60u^8vw \\
                     & + 60u^{6}v^3 + 9u^4 w^2  + 30 u^2 v^2 w + 5v^4,   \\
         \rho_{30}   &= -9x^2 - 12u^{9}vx  -6u^{5}wx  + 9u^{14}vw -10u^{12}v^3  \\
                     &  -3u^{10}w^2  + 30u^8 v^2w   -35u^{6}v^{4} + 6u^{4}vw^{2}  \\
                     & - 10u^2v^3w  -4v^5 -2w^{3}, 
     \end{align*}
where  
  \begin{align*}
      u & = t_{8}, \\
      v & =  2\gamma_{6}+ \gamma_{3}^{2}-u\gamma_{5}+\gamma_{4}(-t^{2}+u^{2})
            -u^{3}\gamma_{3} + t^{6} - t^{4}u^{2}   \\
        &+ t^{3}u^{3} + t^{2}u^{4} - tu^{5}, \\
      w & = \gamma_{10} + u\gamma_{9} - u^{3}c_{7}  -u\gamma_{4}\gamma_{5} 
            +2u^{2} \gamma_{4}^{2} -2u^{2}\gamma_{3}\gamma_{5}  \\
        & + \gamma_{3}\gamma_{4}(-6tu^{2} + 2u^{3}) +  \gamma_{3}^{2}(2t^{2}u^{2} + 2tu^{3} - 2u^{4})  \\
        &  + \gamma_{6}(-5t^{2}u^{2} + 5tu^{3})
         + \gamma_{5}(t^{4}u + 3t^{3}u^{2} + t^{2}u^{3})  \medskip \\
        &  +\gamma_{4}(6t^{4}u^{2} - 3t^{3}u^{3} -2t^{2}u^{4} -tu^{5} + u^{6}) \medskip \\ 
        & +\gamma_{3}(-6t^{5}u^{2} - 2t^{4}u^{3} + 4t^{3}u^{4} + 6t^{2}u^{5} - 4tu^{6} + u^{7})  \medskip \\
        & + 4t^{7}u^{3} -6t^{5}u^{5} + 2t^{4}u^{6}  + t^{3}u^{7} - t^{2}u^{8}, \\
     x &  =  \gamma_{15}-20 \gamma_{3} \gamma_{6}^2 +3  \gamma_{3}^{2}\gamma_{9}  -23 \gamma_{3}^3 \gamma_{6} 
            -6 \gamma_{3}^{5}  +4 \gamma_{6}\gamma_{9}  \\
            &+3 u \gamma_{4}\gamma_{10} - u \gamma_{5}\gamma_{9}  
           -3 u \gamma_{3}^{2} \gamma_{4}^{2} +3 uc_{7}\gamma_{3}\gamma_{4}  -6 u \gamma_{4}^2\gamma_{6}  \\
           & + (-3t + 2 u ) \gamma_{3}^{3} \gamma_{5}    +(-4 t  + 4u) \gamma_{3}\gamma_{5}\gamma_{6}   \\
           &   +(-t^2 - u^2) \gamma_{4}\gamma_{9}  +(t^2 + tu -u^2 )c_{7}   \gamma_{3}^{2}  \\
           & +(9t^2+ 12 tu + 5u^2) \gamma_{3}\gamma_{4}\gamma_{6}  
 \end{align*} 
 \begin{align*} 
  & +  (5t^2 +6tu + 2 u^2)  \gamma_{3}^{3} \gamma_{4}  + (3 t^2 +4tu + u^2)  c_{7} \gamma_{6}        \\
  & +(-6t^3 -2t^2u - 6tu^2 + 5u^3) \gamma_{3}^{4}    - u^3 \gamma_{3}\gamma_{9} \\
  & + (3t^2 u + u^3 ) \gamma_{4}^{3}   + (2 t^2u + 3tu^2) c_{7} \gamma_{5}  \\
  & +(-45t^3 + 10t^2 u -40 tu^2) \gamma_{6}^{2} \\
  &+ (t^3 -2t^2 u + tu^2 - u^3 )\gamma_{3} \gamma_{4} \gamma_{5}  \\
  & +(-33 t^3 + t^2u - 31tu^2 + 13u^3) \gamma_{3}^{2} \gamma_{6} \\
  & +( -2 t^4 - 4 t^3u -3t u^3 + 3u^4) c_{7} \gamma_{4}    \\
  & +(-9 t^4 -6t^3u -18t^2u^2 + 5tu^3 -3u^4)\gamma_{5} \gamma_{6}  \\
  & +(-3t^4 -3t^3 u -7t^2 u^2 + 5tu^3 -4u^4) \gamma_{3}^{2}\gamma_{5} \\
  & +(-t^4 -6 t^3u  -t^2u^2 -3tu^3)\gamma_{3} \gamma_{4}^2    \\
  & + (-3t^4 u -6 t^3 u^2 + 3t^2 u^3  +15 t u^4  )\gamma_{10} \\  
  & + (-3t^4 u +  t^3 u^2 + 5 t^2 u^3 + 10tu^4  -u^5 )c_{7} \gamma_{3}   \\
  &  +(15 t^5 -2t^4u + 3t^3 u^2 + 14t^2 u^3 -16tu^4 + 3u^5) \gamma_{3}^2 \gamma_{4}  \\
  & + (39t^5 - 13t^4u + 8t^3 u^2 + 35 t^2 u^3 -31tu^4 - 3u^5)\gamma_{4} \gamma_{6}     \\
  & +(t^6 -t^4u^2 -{t}^{3}{u}^{3} - t^2 u^4 - tu^5 -u^6)\gamma_{9}    \\
  & +(- 13t^6 + 12t^5 u +5t^4 u^2 -56t^3 u^3 + 8t^2 u^4 +21tu^5  \\
  & + 2u^6) \gamma_{3} \gamma_{6}       \\
  & + ( 6 t^6 + 3t^5 u +2t^4 u^2 +  7 t^3u^3  + t^2 u^4 -8tu^5  + 3u^6)\gamma_{4}  \gamma_{5}    \\
  & +(-8 t^6 + 6t^5u + 2t^4u^2 -22t^3 u^3 + 6t^2u^4 +8tu^5  \\
  & -2u^6) \gamma_{3}^{3}   \\
  & + (-6t^7 + t^6 u -7t^4 u^3  + 5 t^3 u^4  + 3t^2 u^5 + 3tu^6  - 63 u^7)  \gamma_{4}^{2}   \\
  & + (-t^7 +2t^6 u + t^5 u^2 -11 t^4 u^3 + 6t^3 u^4 + 5t^2 u^5 +6tu^6  \\
  & + 39u^7)\gamma_{3}\gamma_{5}   \\
  & + (2t^8 + 6t^7 u +3t^6u^2 -4t^5 u^3  - 15 t^4 u^4  + 6t^3 u^5  \\
  & + 3t^2 u^6  -40tu^7 + 59 u^8)c_{7}   \\
  & + (3t^8 +  t^6 u^2 +11t^5 u^3 + 14t^4 u^4 -20t^3 u^5 -4 t^2 u^6  \\
  & +118tu^7 + 3u^8 )\gamma_{3} \gamma_{4}    \\  
  & + (- 48 t^9 + 3 t^8u  -41 t^7 u^2  + 18 t^6 u^3 + 16 t^5 u^4 \\
  & -13t^4 u^5  - 67t^3 u^6 + 125t^2 u^7 - 15tu^8  -291u^9)\gamma_{6}    \\
  & +(-18t^9 -3t^8 u -16t^7 u^2 + 10t^6u^3 -4t^5u^4 -8t^4 u^5  \\
  & -16t^3 u^6  -23t^2 u^7 -10tu^8  -115u^9)  \gamma_{3}^{2}  \\
  & +(-6t^{10} -3t^9 u -9t^8 u^2 + 5t^7 u^3 -5 t^6 u^4  -14t^4 u^6  \\
  & -52t^3 u^7 + 6t^2 u^8 -60tu^9 +117 u^{10})\gamma_{5}   \\
  & + (18t^{11} -3t^{10} u+ 5t^9 u^2 + 11t^8 u^3 - 28t^7 u^4 + 8t^6 u^5  \\
  & + 20t^5 u^6  -64t^4 u^7 -15t^3 u^8  + 54t^2 u^9 + 178t u^{10}  \\
  & - 177u^{11}) \gamma_{4} \\
  & +(-2t^{12} +6t^{11}u + 2t^{10}u^2 -20 t^9 u^3 +11t^8 u^4  \\
  & + 22t^7 u^5 -8t^6 u^6 + 83t^5 u^7  + 15t^4 u^8 + 5t^3 u^9  \\
  & -116 t^2 u^{10}  
    + tu^{11} + 117 u^{12}) \gamma_{3}        \\
 &  -12t^{15} - t^{14}u -10t^{13}u^2 + 6 t^{12} u^3 + 7t^{11}u^4   -13t^{10}u^5  \\
 & -31t^9 u^6  +9 t^8 u^7 -t^7 u^8 -118t^6 u^9    -18t^5 u^{10}  \\
 &+ 131t^4 u^{11}  -6t^3 u^{12}   - 233t^2u^{13}  + 175tu^{14} -58 u^{15}.        
  \end{align*} 
\end{thm}

\section{Concluding remarks} 
The flag manifold $G/T$ also appears as $G_{\C}/B$, where $G_{\C}$ denotes the complexification of $G$ and 
$B$ a Borel subgroup of $G_{\C}$ containing $T_{\C}$, the complexification of $T$.  The subgroup $B$ acts on 
$G_{\C}/B$ by the  left translation.     Each element $w \in W$ defines an element $e_{w} := wB$ in 
$G_{\C}/B$,  and the $B$-orbit $X_{w}^{\circ} := B \cdot e_{w}$ of $e_{w}$ under this action  
is isomorphic to an affine space,   and  is called a {\it Schubert cell} corresponding to $w$.  
Then it was shown by Chevalley \cite{Che58}  
that the manifold $G_{\C}/B$ admits a cellular  decomposition 
\begin{equation}  \label{eqn:cell_decomposition} 
          G_{\C}/B = \coprod_{w \in W} X_{w}^{\circ}.    
\end{equation}   
The {\it Schubert variety} $X_{w}$ corresponding to $w$ is defined to be the  closure of  the  Schubert cell $X_{w}^{\circ}$, 
and it determines a {\it Schubert class} in $H^{*}(G_{\C}/B;\Z)$.  
It follows from (\ref{eqn:cell_decomposition}) that the set of Schubert classes form a $\Z$-basis for 
$H^*(G_{\C}/B;\Z)$ (for details, see \cite{BGG73}, \cite{Che58},  \cite{Ful97}).  This is the {\it Schubert presentation} of 
the cohomology ring  of $G/T$.  It is natural to ask the connection between the Borel and the Schubert presentations. 
In \cite{BGG73}, Bernstein, Gelfand and Gelfand introduced the {\it divided difference operators} and gave 
a  general answer to the above problem.   In our situation, the generators $t_{i} \; (1 \leq i \leq 8)$, $t$ and 
$\gamma_{i} \; (i = 3, 4, 5,6, 9,10,15)$ in Theorem \ref{thm:E_8/T} can be  expressed as $\Z$-linear combinations 
of Schubert classes indexed by   $W(E_{8})$. 
The explicit forms will be given in our forthcoming paper (\cite{Kaji-Nak2}) (see also \cite{DZ08}
for the Schubert presentation of $H^*(E_{8}/T;\Z)$). 


\end{document}